\documentclass[12pt]{amsart}\usepackage[v2]{xy}
\def\b#1{{\boldsymbol{#1}}}
\def\btau{\bar{\tau}}

\def\dpar{\partial}
\def\lra{\to}

\newtheorem{thm}[equation]{Theorem}
\newtheorem{cor}[equation]{Corollary}
\newtheorem{lem}[equation]{Lemma}
\newtheorem{prop}[equation]{Proposition}
\theoremstyle{definition}
\newtheorem{defn}[equation]{Definition}
\newtheorem{conj}[equation]{Conjecture}
\theoremstyle{remark}
\newtheorem{rem}[equation]{Remark}

\numberwithin{equation}{section}

\def\koszul#1#2{{\epsilon_{#1,#2}}\;}
\def\Lie#1{{\sf Lie}(#1)}
\def\Z{\mathbb{Z}}

\begin{document}
\title{Homotopy Gerstenhaber Structures and Vertex Algebras}
\author{I G\'alvez \and V Gorbounov\and A Tonks}
\address{London Metropolitan University}\email{i.galvezicarrillo@londonmet.ac.uk}
\address{University of Aberdeen}\email{vgorb@maths.abdn.ac.uk}
\address{London Metropolitan University}\email{a.tonks@londonmet.ac.uk}

\thanks{The first and third authors were partially supported by
 MEC/FEDER grants MTM2004-03629 and MTM2007-63277. The second author was partially supported by the NSF (contract 0106574)}
\begin{abstract}
We provide a simple construction of a $G_\infty$-algebra structure 
on an important class of vertex algebras $V$,
which lifts the Gerstenhaber algebra structure on BRST cohomology of $V$
introduced by Lian and Zuckerman.
We outline two applications to algebraic topology: the construction of
a sheaf of $G_\infty$ algebras on a Calabi--Yau manifold $M$,
extending the operations of multiplication and bracket of functions
and vector fields on $M$, and of a Lie$_\infty$ structure related to
the bracket of Courant~\cite{courant}. 
\end{abstract} 
\maketitle

\section*{Introduction}

The idea that algebraic structure present on the (co)homology of certain objects 
may be the shadow of some richer `strong homotopy algebra' structure is a recurring theme in several areas of mathematics. 
The first example was the recognition theorem of Stasheff, which says that a CW-complex is a loop space if and only if it has an $A_\infty$ structure. 
More recently Deligne conjectured that the Gerstenhaber algebra structure 
on the Hochschild cohomology of an algebra could be lifted to a $G_\infty$ structure on the Hochschild complex itself, and this was demonstrated in a number
of papers~\cite{bergerfresse,ks00,mccluresmith,hgavoronov}.

In this paper we prove a similar conjecture inspired by mathematical physics.
The fact that the BRST cohomology of a topological conformal field theory is
a Gerstenhaber algebra was observed by Lian and Zuckerman \cite{lianzuck} in 1993.
They defined product and bracket operations on the corresponding chain complex which satisfy the Gerstenhaber axioms up to homotopy and suggested the existence of `higher homotopies'. In \cite[Conjecture 2.3]{kivozu} the question was refined to what is now known as the Lian--Zuckerman conjecture: the product and bracket operations on a topological vertex operator algebra can be extended to the structure of a $G_\infty$ algebra.

To motivate the conjecture further, recall that 
$\Z_{\geq 0}$-graded  vertex algebras carry a
commutative algebra structure on the weight zero part, 
see for example~\cite[p.\ 623]{chiral2}. In
Lemma~\ref{below} we see that if a vertex algebra has conformal weight
bounded below, then the part of lowest weight in fact carries a Gerstenhaber
algebra structure. It is therefore natural to ask what structure can
be found on the components of higher weight.

We prove the Lian--Zuckerman conjecture in Theorem~\ref{maintheorem} below for a wide class of vertex algebras, including positive definite lattice vertex algebras and Kac--Moody algebras. Our motivating application is to the vertex algebra structure on a manifold given by the chiral de Rham complex, which was introduced in \cite{msv} and whose relation to string theory is outlined in \cite{kap,witten}. 
Using methods analogous to those developed by Stasheff {\em et al.}~\cite{shlie} we are able to construct a canonical $G_\infty$ algebra structure
in these cases. 
We point out in Remark~\ref{obstruct} that, in all cases, the existence of such a structure is equivalent to the vanishing of certain maps $\Gamma_{r+1}$ in cohomology.

A related problem was investigated by Huang and Zhao \cite{huzh}, where it was shown that an appropriate topological completion of a topological vertex operator
algebra has a structure of a genus-zero topological conformal field theory (TCFT). Together with the result that a TCFT has a natural $G_\infty$ structure claimed in~\cite{kivozu} this would settle the Lian--Zuckerman conjecture, but the $G_\infty$ algebras as defined there do not exist, due to an error in the highly influential preprint of Getzler and Jones~\cite{getzler-jones} pointed out by Tamarkin. For more details, and a proposed correction to~\cite{kivozu}, we refer the reader to Voronov's discussion in~\cite[Section 4]{hgavoronov}.

Tamarkin and Tsygan~\cite{tt00} proposed an explicit definition of $G_\infty$ algebras based on the minimal model of the Gerstenhaber operad given by Koszul duality, and this is the definition we take here. They also gave a very general definition of homotopy Batalin--Vilkovisky algebra, to which we will return in a later paper.

We will say a few words about the applications of our construction to 
geometry and topology.
The construction of the chiral de Rham complex due to 
Malikov, Schechtman and Vaintrob~\cite{msv}
provided a new algebraic structure, namely a topological vertex algebra, 
which is associated to any  manifold and extends the operations of 
multiplication of functions and of the Lie bracket of vector fields.
It is therefore natural to consider the conjecture of Lian--Zuckerman in this mathematical context. This is another task we address in this paper: 
to a manifold we can associate a sheaf of $G_\infty$-algebras derived from the chiral de Rham complex. 

The connection between algebras up-to-homotopy and structures related to vertex algebras has also been observed in a slightly different guise by Roytenberg and Weinstein, who show in~\cite{rw-courant} that a Courant algebroid structure leads naturally to that of a homotopy Lie algebra. Our result is a generalisation of this observation, since according to Bressler~\cite{bressler-p1} a Courant algebroid is a quasi-classical limit of a vertex algebroid.

A short outline of the paper is as follows. We first recall basic conventions for graded and super algebra and the definition of Gerstenhaber and $G_\infty$ algebra. We develop a notion of partial $G_\infty$ structure and an abstract algebraic situation of which the Lian--Zuckerman situation is the motivating example. We prove the existence of a $G_\infty$ structure in this general context, and give an explicit formula for this structure. In the next section we turn to vertex algebras, recalling the necessary definitions and applying our results to prove the Lian--Zuckerman conjecture for vertex algebras with non-negative conformal weights. We end by indicating some important examples of such vertex algebras which arise from the chiral de Rham complex and chiral polyvector fields.

We would like to acknowledge helpful conversations with Jim Stasheff and Dmitri Tamarkin. 
The authors visited the Max-Planck-Institut f\"ur Mathematik and
the Universitat de Barcelona while working on this project and are grateful to
these institutions for excellent working conditions.

\section{Graded (or super) algebra} 

We begin by recalling some very basic notions of graded algebra and superalgebra.

All vector spaces are over a field of characteristic zero. If
$W=\bigoplus_{n\in\Z} W_n$ is a graded vector space then $W[k]$ denotes the $k$th
desuspension, defined by $$(W[k])_n=W_{n+k}.$$
The notation $s^{-1}W$ instead of $W[1]$ is common, especially in algebraic topology.

A graded linear map
$f:V\to W$ of degree $k$ between graded vector spaces 
is a linear map $f:V\to W[k]$ with $f(V_n)\subseteq W[k]_n$.
The degree of a homogeneous element or graded linear
map $a$ is denoted by $|a|$, and we write $\koszul a b$ for the sign $(-1)^{|a|\,|b|}$.
 The natural symmetry isomorphism of graded vector spaces is given by
\begin{align*}
V\otimes W&\to \;\;\;W\otimes V\\
a\otimes b&\mapsto \koszul a b b\otimes a
\end{align*}
for homogeneous elements $a,b$ in $V,W$. 

If the vector space $W$ is not $\Z$- but $\Z_2$-graded then $W$ is termed a superspace,
and we denote by $W_{ev}$ and $W_{odd}$ the subspaces of even and odd vectors. The degree $|a|$ of an element is often written $p(a)$ and termed parity in this situation.

The {\em commutator} with respect to a binary operation $\cdot$
on a graded (or super) vector space $W$ is defined by
\[
[a,b]\;\;=\;\;a\cdot b-\koszul a b b\cdot a\,.
\]
The operation $\cdot$ is
\begin{align*}
\text{\em commutative\;\;}&\text{if\;\;}[a,b]=0\\
\text{\em skew-symmetric\;\;}&\text{if\;\;}a\cdot b+\koszul a b b\cdot a=0
\end{align*}
for homogeneous elements $a,b$ in $W$. 
A graded {\em Lie algebra} is a graded vector space $W$ with
a 
degree zero skew-symmetric bilinear operation $[\:,\,]\colon W\otimes W\to W$ 
which satisfies the graded Jacobi identity
\[
[[a,b],c]+\koszul a b\koszul a c[[b,c],a]+\koszul a c\koszul b c [[c,a],b]=0
\]
for homogeneous elements $a,b,c$ in $W$.
A graded linear map $f:W\to W$ is a {\em derivation} with respect to a
binary operation $\cdot$ on $W$ if it satisfies
\begin{align*}
f(a\cdot b)&=f(a)\cdot b\,+\,\koszul f a a\cdot f(b)\,.
\end{align*}
for homogeneous elements $a,b\in W$.

A {\em Gerstenhaber algebra} is a graded vector space $W$ with
bilinear operations $\cdot$ of degree zero and $[\;,\,]$ of degree $-1$ such that
\begin{enumerate}
\item
$(W,\cdot)$ is a graded commutative associative algebra
\item
$(W[1],[\;,\,])$ is a graded Lie algebra
\item
$[\,a,{-}]:W\to W$ is a derivation with respect to the bilinear operation $\cdot$, for each
homogeneous $a\in W$.
\end{enumerate}

We note the following folklore result 
for the commutator with respect to composition.

\begin{lem}\label{deriv}
Suppose that $f,g:W\to W$ are derivations with respect
  to a bilinear operation $\cdot$ on a graded vector space $W$.
Then 
\begin{enumerate}\item
the graded commutator $[f,g]:W\to W$ is a derivation,
\item if $f$ has odd degree then $f\circ f:W\to W$ is a derivation.
\end{enumerate}
\end{lem}
\begin{proof}
(1) Let $a,b$ be homogeneous elements of $W$. Then
\begin{eqnarray*}
[f,g](a\!\cdot\! b)&\!\!\!=\!\!\!&
\;\;f(g(a\!\cdot\! b))\;\;-\;\;
\koszul fg g(f(a\!\cdot\! b))
\\&\!\!\!=\!\!\!&
f\left(ga\!\cdot\! b+\koszul ga a\!\cdot\! gb\right)
-\koszul fg
g\left(fa\!\cdot\! b+\koszul fa a\!\cdot\! fb\right)
\\&\!\!\!=\!\!\!&
fga\!\cdot\! b + \koszul f{ga}
ga\!\cdot\! fb + \koszul ga
fa\!\cdot\! gb + \koszul ga\koszul fa
a\!\cdot\! fgb
\\
-\koszul f g
\!\!\!\!\!\!
&\!\!\!\bigl(\!\!\!\!\!\!\!\!&
gfa\!\cdot\! b + \koszul g{fa}
fa\!\cdot\! gb + \koszul fa
ga\!\cdot\! fb + \koszul fa\koszul ga
a\!\cdot\! gfb
\;\bigr)
\\&\!\!\!=\!\!\!&
[f,g]a\!\cdot\! b 
\quad+\quad 0\quad+\quad 0\quad+\quad 
\koszul{[f,g]}a\,a\!\cdot\! [f,g]b
\end{eqnarray*}
since $\koszul f{ga}=\koszul fg\koszul fa$ etc.

(2) This follows by a similar explicit calculation. Alternatively, since we are working over a field of characteristic zero, we note that if $f$ is odd then $f\circ f=\frac12[f,f]$ is a derivation by part (1).
\end{proof}

\section{$G_\infty$-algebras}\label{Ginf}

We recall the definition of a $G_\infty$ algebra from~\cite{tt00}.
Let $\Lie A$ be the free graded Lie algebra generated by a graded vector space $A$, 
\[
\Lie A
\;=\;\bigoplus_p L^pA,\qquad
L^pA\;=\;
[[[\ldots\![A,A],\ldots A],A],A]
\]
where $L^pA$ is spanned by all $(p-1)$-fold commutators of elements of $A$.
Then consider the free graded commutative algebra on the suspension  $(\Lie A)[-1]$, which we can write as
\[
GA\;\;=\;\;\bigoplus_t\bigwedge{\!\!}^t\,\Lie A\,[-t].
\]
The Lie bracket on $\Lie A$ extends to a degree $-1$ bilinear operation 
$$[\:,\,]\colon GA\otimes GA\to GA$$
which is a derivation with respect to the $\wedge$-product on $GA$.

There are several formulations of the notion of $G_\infty$ structure
in the literature. We follow that of Tamarkin
and Tsygan in \cite[Section 1]{tt00}. 

\begin{defn}
A {\em $G_\infty$ structure} on a graded
vector space $V$ is a square zero degree one linear map
\[
\gamma:G(V[1]^*)\to G(V[1]^*)
\]
which is a derivation with respect to both $\wedge$ and $[\;,\,]$.
\end{defn}

\begin{rem}
We have implicitly assumed $V$ is finite dimensional; in particular 
the grading on $V$ is bounded so
that the linear dual $A=V[1]^*$ is also graded, with 
$A_n={V_{1-n}}^*$.
This is not really a restriction, but just a side-effect of defining
$G_\infty$ structures on $V$ as derivations on free algebras on
$V^*$. To avoid the finite dimensionality assumption we could instead work with
{\em co}derivations on free {\em co}algebras on $V$; compare for
example \cite[Section 1.2.2]{dtt07}. 
\end{rem}

Consider the components
\[
G^{p_1,p_2,\ldots,p_t}A
\;=\;L^{p_1}A\wedge L^{p_2}A \wedge \cdots\wedge L^{p_t}A \,[-t]
\;\subseteq\;GA
.
\]
In fact it is slightly more manageable to consider
\[
G_m
A\;\;=\;\;\bigoplus_{
\substack{
p_1+\cdots+ p_t=m}
} G^{p_1,p_2,\ldots,p_t}A\]
so that $G_1A=A[-1]$ and $G_2A=[A,A][-1]\oplus A\wedge
A\,[-2]$, etc.
An element of $G_mA$ is said to have {\em length} $m$ in $GA$. 
\begin{lem}
There are bijections between the set of degree 1 linear maps $\gamma\colon GA\to GA$
which are derivations with respect to both $\wedge$ and $[\;,\,]$, 
the set of degree 1 linear maps $\gamma_1\colon G_1A\to GA$, and the set of sequences of degree 1 linear
maps $\gamma_1^{m+1}\colon G_1A\to G_{m+1}A$, $m\ge0$.
\end{lem}
\begin{proof}Any such derivation $\gamma$ is determined
by its restriction $\gamma_1=\gamma|_{G_1A}$, 
and $\gamma_1$ is the sum of its components
$\gamma_1^{m+1}\colon G_1A\to G_{m+1}A$.
\end{proof}
Derivations behave well with respect to the filtration by
length. For fixed $m\geq0$ the derivation laws allow us to extend the linear map
\[
\gamma_1^{m+1}\colon G_1A\to G_{m+1}A
\]
in the lemma to a family of linear maps $\gamma_i^{m+i}$, for $i\ge1$, 
\[
\gamma_i^{m+i}\colon G_iA\to G_{m+i}A.
\]
The $G_\infty$-algebra condition that a derivation
$\gamma:GA\to GA$ has square zero is clearly equivalent to the
following collection of quadratic relations on these components,
{\renewcommand{\theequation}{{\sf R}$_i^k$}
\begin{eqnarray}\sum_{j=i}^k \gamma_j^k\gamma_i^j&=&0
\;:\;  G_iA\to G_kA
.
\end{eqnarray}
}

For our purposes it will be useful to introduce 
a notion of a {\em partial} $G_\infty$ structure.

\begin{defn} Let $V$ and $A$ be as above and let $1\leq r\leq \infty$. 
A {\em $G_r$ structure} $\gamma^{\le r}$ 
on $V$ is a sequence of degree 1 linear maps 
$\gamma_1^k\colon G_1A\to G_kA$ 
for $k-1< r$ which, together with their extensions
$\gamma_i^j:G_iA\to G_jA$, satisfy
the relations {\sf R}$_1^k$ for $k-1< r$.  
\end{defn} 

The definition requires the relations ${\sf R}_i^k$ to hold only when
$i=1$; to justify the name we note that corresponding `higher' relations hold automatically:

\begin{lem}
A $G_r$-structure $\gamma^{\le r}$  on $V$ satisfies the relation {\sf
  R}$_i^k$ whenever $k-i<r$.
\end{lem}
\begin{proof}
Extend each map $\gamma_1^j$ to a degree one derivation $GA\to GA$,
$j\le r$, and let $f$ be the sum of these derivations.
Now consider the map $g=f^2:GA\to GA$, 
which by Lemma \ref{deriv}(2) is a also derivation. 
If $1\le k\le r$ the relation {\sf R}$_1^k$
says that the component $g_1^k:G_1A\to G_kA$ is zero.
Hence its extensions  
$g_i^{i+k-1}:G_iA\to G_{i+k-1}A$ are zero and the relations 
{\sf R}$_i^{i+k-1}$ hold also.
\end{proof} 
For $r=\infty$ the new notion of $G_r$ structure is consistent with the previous definition:
\begin{cor} 
A $G_\infty$ structure on a graded vector space $V$ is specified by a family of degree 1 linear maps 
$$\gamma_1^k:G_1(V[1]^*)\to G_k(V[1]^*),\qquad k\geq1,$$ 
which (together with their extensions $\gamma_j^k$) satisfy the relations ${\sf R}_1^k$, $k\geq1$.
\end{cor}

In the case $r=1$ a $G_r$-structure is just a degree one map
$$d_1:=\gamma_1^1:G_1A\to G_1A$$ 
which has square zero, 
and the Lemma says the extensions 
$$d_m:=\gamma_m^m:G_mA\to G_mA$$ 
also have square zero.

For $k\ge2$ we can express the relation {\sf R}$_1^k$ in the form
\begin{align}\label{nullhtpy}
\gamma_1^kd_1
\;+\;
d_k\gamma_1^k
&=\Gamma_k
\end{align}
where
\begin{align*}
\Gamma_k&=-
\sum_{j=2}^{k-1}
\gamma_{j}^{k}\gamma_{1}^{j}
\;:\;  G_1A\to G_kA
.\end{align*}
Note that,
for any $G_r$ structure $\gamma^{\leq r}$, 
the linear map
$\Gamma_{r+1}$ is defined. The following
result, that it is always a chain map $(G_1A,d_1)\to(G_{r+1},d_{r+1})$,
is a translation of~\cite[Lemma 6]{shlie} from the Lie$_\infty$ to the $G_\infty$ context.
\begin{prop}\label{stashtrick}
Suppose $\gamma^{\le r}$ is a $G_r$ structure. Then 
the linear map $\Gamma_{r+1}$ satisfies
\begin{eqnarray*}
d_{r+1} \Gamma_{r+1}
&=&
\Gamma_{r+1} d_1.
\end{eqnarray*}
\end{prop}
\begin{proof}
We have
\begin{eqnarray*}
d_{r+1}\Gamma_{r+1}= -\!\!\!
\sum_{1<i<r+1}\!\!\!\!\gamma_{r+1}^{r+1}\gamma_i^{r+1}\gamma_1^i\!\!&=&
\!\!\!\!\!\!\!\sum_{1<i\le j<r+1}\!\!\!\!\!\!\!\gamma_{j}^{r+1}\gamma_i^{j}\gamma_1^i
\;\text{ by the relations \sf R}_i^{r+1},\\
\Gamma_{r+1}d_1=     -\!\!\!
\sum_{1<j<r+1}\!\!\!\!\gamma_{j}^{r+1}\gamma_1^{j}\gamma_1^1\;\;&=&
\!\!\!\!\!\!\!\sum_{1<i\le j<r+1}\!\!\!\!\!\!\!\gamma_{j}^{r+1}\gamma_i^{j}\gamma_1^i
\;\text{ by the relations \sf R}_1^{j}.
\end{eqnarray*}
\end{proof}

\begin{rem}\label{obstruct}
Now equation (\ref{nullhtpy}) can be read as saying that the obstruction to
extending a $G_r$ to a
$G_{r+1}$ structure is the existence of 
\begin{itemize}
\item a null-homotopy $\gamma_1^{r+1}$ for the chain map $\Gamma_{r+1}$, 
\item or an element $\gamma_1^{r+1}\in{\mathrm{Hom}}(G_1A,G_{r+1}A)$ whose coboundary is the cocycle $\Gamma_{r+1}$.
\end{itemize}
\end{rem}

\section{Main theorem}

We now set up an abstract situation where we can apply the above results to define explicit $G_\infty$ structures. In the following section we show how this applies to the Lian--Zuckerman conjecture. 

Consider a bigraded vector space 
\begin{equation}\label{bigraded}
V=\bigoplus_{s\geq0} V^s
=\bigoplus_{s\geq0} \left(
\bigoplus_{N'_s\leq n\leq N_s} V_n^s\right)
\,.
\end{equation}
We refer to the two gradings as the {\em fermionic} degree, written $|v|=n$, and the
{\em conformal} degree, written $\|v\|=s$. As in section \ref{Ginf} we define
$$GA=\bigwedge{\!\!}^\bullet((\Lie A)[-1]),\qquad A=V[1]^*\,,$$
where the (de)suspensions are with respect to the fermionic degree.

Let us assume for a moment that the conformal degree on $V$ is also bounded
above. Then 
$A$ is isomorphic to the (finite) direct sum of $A^s:=V^s[1]^*$, and
this conformal grading extends to $GA$ by 
$$\|v\wedge w\|=\|[v,w]\|=\|v\|+\|w\|\,.$$

All linear maps we consider in this section will be of degree zero with respect to the conformal degree; in
particular only the fermionic degree will contribute to the signs in
derivation formulas. In the interests of brevity, and consistency with
section \ref{Ginf}, 
the word ``degree'' will mean fermionic
degree unless otherwise stated.

\begin{thm}\label{main}
Suppose $V$ is a bigraded vector space as above together
with linear maps
$$
m_1\colon V\to V,\qquad m_2,\;m_{1,1}\colon V\otimes V\to V
$$
whose duals define a $G_2$ structure on $V$, and a square zero
linear map
$$
h\colon V\to V
$$
such that, for all elements $v$ of conformal degree $s$ in $V$,
$$m_1hv+hm_1v=sv.$$
Then any extension of the $G_2$ structure to a
$G_\infty$ structure on $V^0$ has an extension to a 
$G_\infty$ structure on $V$.
\end{thm}
\begin{proof}
Extending the duals of the maps $m_1$ and $h$ to derivations on
$GA$ we obtain square zero linear maps $d$ and $\sigma$, of degrees 1 and $-1$
respectively, satisfying
\begin{align}\label{dsig}
d\sigma a+\sigma da&=sa
\end{align}
if $\|a\|=s$. The map $\sigma$ may be thought of as a chain homotopy.

Let $s>0$ and suppose inductively we have: a $G_r$ structure on $V$, a
$G_{r+1}$ structure on the subspace of elements of conformal
degree $<s$, and linear maps $\gamma_1^{r+1}\colon G_1A\to G_{r+1}A$
defined on elements of conformal degree $s$ and fermionic degree $\leq n$
which satisfy
\begin{eqnarray}\label{rel}
\gamma_1^{r+1}da+
d\gamma_1^{r+1}a
&=&\Gamma_{r+1}a
\end{eqnarray}
if $\|a\|=s$ and $|a|<n$.
Now if $\|a\|=s$ and $|a|=n+1$, define
\begin{align}\label{gdef}
\gamma_1^{r+1}a
&=
\frac1s(
\Gamma_{r+1}\sigma a-d\gamma_1^{r+1}\sigma a
)
\end{align}
We must show the relation (\ref{rel}) holds if $\|a\|=s$ and $|a|=n$.
We have
\begin{eqnarray*}
\gamma^{r+1}_1da&=&
\frac1s(
\Gamma_{r+1}-d\gamma_1^{r+1}
)\sigma(da)
\\&=&
(
\Gamma_{r+1}-d\gamma_1^{r+1}
)(1-\frac1sd\sigma)(a)
\end{eqnarray*}
by equations $(\ref{dsig})$ and $(\ref{gdef})$, so that
\begin{eqnarray*}
\gamma_1^{r+1}da
\;-\;\Gamma_{r+1}a\;+\;
d\gamma_1^{r+1}a
&=&
-\frac1s(
\Gamma_{r+1}-d\gamma_1^{r+1}
)d\sigma a
\\&=&
-\frac1sd(
\Gamma_{r+1}-\gamma_1^{r+1}d
)\sigma a
\\&=&
-\frac1sd(
d\gamma_1^{r+1}
)\sigma a
\;\;=\;\;0
\end{eqnarray*}
using Proposition \ref{stashtrick}, the inductive hypothesis and $d^2=0$.
\end{proof}

\section{Vertex algebras and the structure of BRST cohomology}

There are a number of expositions of vertex algebra theory; we include the following for completeness from \cite{chiral3}, which is also the source of our examples to which we return in the final section.

A {\it $\Z_{\geq 0}$-graded 
vertex superalgebra} (over a field $k$) is a 
$\Z_{\geq 0}$-graded superspace $V=\bigoplus_{i\geq 0}V^i$, where the component in each conformal degree $i$ has a parity decomposition $V^i=V^i_{ev}\oplus V^i_{odd}$,  equipped 
with a distinguished {\it vacuum vector}  $\b1\in V^0_{ev}$
and a family of bilinear operations 
$$
_{(n)}\colon V\times V\to V\,,
$$
for $n\in\Z$, such that 
$$
p(a_{(n)}b)=p(a)+p(b)\,,\qquad {V^i}{}_{(n)}{}V^j\subset V^{i+j-n-1}\,.
$$
The following properties must hold:
$$
\b1_{(n)}a=\delta_{n,-1}a\,,\qquad
a_{(-1)}\b1=a\,,\qquad
a_{(n)}\b1=0\ \text{for }n\geq 0 
$$
and 
\begin{align*}
&\sum_{j=0}^\infty\binom{m}{j}(a_{(n+j)}b)_{(m+l-j)}c=
\\&=\sum_{j=0}^\infty(-1)^j\binom{n}{j} 
\bigl\{a_{(m+n-j)}(b_{(l+j)}c)-(-1)^{n+p(a)p(b)}b_{(n+l-j)}(a_{(m+j)}c)\bigr\}
\end{align*}
for all $m,n,l\in \Z$.

A {\em topological} vertex algebra is a
vertex superalgebra $V$ with certain extra structure. The full mathematical definition may be found in \cite[section 5.9]{kac}. For our purposes it is enough that there are distinguished elements $L\in V^2_{ev}$, $G\in V^2_{odd}$, $Q\in V^1_{odd}$,  $J\in V_{ev}^1$ satisfying
$$
{Q_{(0)}}^2=0={G_{(1)}}^2,\qquad\qquad
[Q_{(0)},G_{(i)}]=L_{(i)}
\,.
$$ 
Moreover $L_{(1)}$ and $J_{(0)}$ are commuting diagonalizable operators whose eigenvalues coincide with the conformal and fermionic gradings. Modulo 2, the fermionic grading gives the parity.
In the terminology  of conformal field theory, the elements $L$, $G$, $Q$ and $J$ are known to physicists as the Virasoro element, the superpartners and the current. They have fermionic gradings 0, $-1$, 1 and 0 respectively.

Lian and Zuckerman in \cite{lianzuck} observed the following

\begin{thm} 
On any topological vertex algebra $V$ the operations 
$$x\,\bullet\,y=x_{(-1)}y,\qquad\qquad
\{x,y\}=(-1)^{p(x)}(G_{(0)}x)_{(0)}y$$ 
are cochain maps with respect to the differential
$d=Q_{(0)}$ and induce a Gerstenhaber algebra structure on the
cohomology  $H^*(V,d)$.
\end{thm}

They posed the following conjecture of lifting the
Gerstenhaber algebra structure on the cohomology to a homotopy algebra
structure on $V$ itself, compare \cite[Conjecture 2.3]{kivozu} and \cite{lianzuck}.

\begin{conj}\label{lzconj}
Let $V$ be a topological vertex algebra. Then the product and bracket
\begin{align*}
x\cdot y&={\textstyle\frac12}(x_{(-1)}y +(-1)^{p(x)p(y)} y_{(-1)}x)
\\
(-1)^{p(x)}[x,y]&={\textstyle\frac12}( (G_{(0)}x)_{(0)}y+(-1)^{p(x)p(y)} (G_{(0)}y)_{(0)}x )
\end{align*}
extend to a $G_\infty$ structure on $V$.
\end{conj}
\noindent This product $x\cdot y$ is by
definition graded commutative and of degree zero. 
The bracket is a graded skew-symmetric operation of fermionic degree
1,
$$
[x,y]={\textstyle\frac12}(\{x,y\}-(-1)^{(p(x)-1)(p(y)-1)}\{y,x\}).
$$
Both operations have conformal weight zero. 

It is a well known result that the $\,{}_{(-1)}\,$ operation on a 
$\Z_{\geq 0}$-graded 
topological vertex
algebra $V$ gives the structure of a
graded commutative algebra on $V^0$, 
see for example~\cite[p.623]{chiral2}. This result can be extended as follows:
\begin{lem}\label{below}
Suppose $V$ is a topological vertex algebra with
conformal weights bounded below, $V=\oplus_{i\geq k}V^i$. 
Then the product $\;\cdot\;$ and the bracket $[\;,\,]$ in
Conjecture~\ref{lzconj} above equip $V^k$ with the structure of a
differential Gerstenhaber algebra.
\end{lem}
\begin{proof}
Lian and Zuckerman show that up to homotopy the operation 
$\bullet$ is commutative and associative and the operation $\{\;,\:\}$ is
skew-symmetric, satisfies the Jacobi
identity, and is a derivation of $\bullet$.
Explicit chain homotopies of conformal weight zero are given for these
laws, in equations (2.14), (2.16), (2.23), (2.25), (2.28)
of~\cite{lianzuck} respectively. 
Applied to elements of conformal weight $k$, all of these homotopies
factor through terms of weight $k-1$; hence they
are trivial and the Gerstenhaber axioms hold on the nose in $V^k$. 
In particular, the operations $\bullet$ and $\{\;,\:\}$  
are respectively commutative and skew-symmetric on $V^k$ and
coincide with
the operations $\cdot$ and $[\;,\:]$ in Conjecture \ref{lzconj}.
\end{proof}

We can now apply Theorem \ref{main} to Conjecture~\ref{lzconj}. A
topological vertex algebra $V$ is bigraded by the eigenvalues of $L_{(1)}$
and $J_{(0)}$, and
$$m_1=Q_{(0)}:V\to V,\qquad
m_2=\;\cdot\;,\;\;m_{1,1}=[\;,\,]:V\otimes V\to V$$
give a $G_2$ structure on $V$ and, if it is $\Z_{\geq 0}$-graded, a Gerstenhaber structure on $V^0$. Also
$$h=G_{(1)}:V\to V$$
is a linear map which satisfies
$
[Q_{(0)},G_{(1)}]=L_{(1)}\,$.
Hence we have: 

\begin{thm}\label{maintheorem}
Any $\Z_{\geq0}$-graded topological vertex algebra, such that for each
conformal weight the fermionic grading is finite, has an explicit
$G_\infty$ structure which extends the Gerstenhaber structure on $V^0$ and reduces to the Lian--Zuckerman structure in cohomology.
\end{thm}

There is a wealth of examples of such vertex algebras, including the
chiral de Rham and polyvector fields which we will discuss in the next section.

\section{Applications}

As we observed in the introduction, the work of~\cite{msv} allows us to consider the Lian--Zuckerman conjecture in the context of algebraic topology.
We cannot attempt to give full details but we advise the reader to consult the papers \cite{chiral2,chiral3}.

In \cite{chiral2} a method of constructing vertex algebras was
introduced, starting with a so-called vertex algebroid $L$ and
applying to it the vertex envelope construction, analogous to the
universal envelope for Lie algebras. Important examples of vertex
algebras arise in this way. We will describe here the vertex
algebroids $L$ whose vertex envelopes are the vertex algebras of chiral de Rham
differential forms and chiral polyvector fields. The construction of
the vertex envelope is fairly straightforward; details may be
found in \cite{chiral3}. 

An {\it extended Lie superalgebroid}  is a quintuple 
${\mathcal T}=(A,T,\Omega,\dpar,\langle,\rangle)$ where 
$A$ is a supercommutative $k$-algebra, 
$T$ is a Lie superalgebroid over $A$, 
$\Omega$ is an $A$-module equipped with a structure of a module over the Lie superalgebra $T$, 
$\partial\colon A\lra\Omega$ is an even $A$-derivation and a morphism of $T$-modules, 
and
$\langle,\rangle\colon T\times\Omega\to A$ is an even $A$-bilinear pairing. 

The following identities must hold, for $a\in A$, $\tau,\nu\in T$, $\omega\in\Omega$:
\begin{align*}
\langle\tau,\dpar a\rangle&=\tau(a)
,
\\
\tau(a\omega)&=\tau(a)\omega+(-1)^{p(\tau)p(a)}a\tau(\omega)
,
\\
(a\tau)(\omega)&=a\tau(\omega)+
(-1)^{p(a)(p(\tau)+p(\omega))}\langle\tau,\omega\rangle\dpar a
,
\\
\tau(\langle\nu,\omega\rangle)&=\langle [\tau,\nu],\omega\rangle+ 
(-1)^{p(\tau)p(\nu)}\langle\nu,\tau(\omega)\rangle
.
\end{align*}

Now a {\it vertex superalgebroid} is a septuple $(A,T,\Omega,\partial,
\gamma,\langle\ ,\ \rangle,c)$ where $A$ is a supercommutative $k$-algebra, 
$T$ is a Lie superalgebroid over $A$, $\Omega$ is an $A$-module equipped 
with an action of the Lie superalgebra $T$, $\dpar\colon A\lra\Omega$ is an even  
derivation commuting with the $T$-action, 
$$
\langle\ ,\ \rangle\colon(T\oplus\Omega)\times (T\oplus\Omega)\lra A
$$
is a supersymmetric even $k$-bilinear pairing equal to zero on 
$\Omega\times\Omega$ and such 
that $
(A,T,\Omega,\dpar,\langle\ ,\ \rangle|_{T\times\Omega})$ 
is an extended Lie superalgebroid, $c\colon T\times T\lra\Omega$ is a 
skew supersymmetric even $k$-bilinear pairing and 
$\gamma\colon A\times T\lra\Omega$ is an even $k$-bilinear map. 
The maps $c$ and $\gamma$, and the pairing $\langle\:,\;\rangle$, satisfy a number of axioms which essentially express their failure to be $A$-linear \cite{chiral3}.

Let $A$ be a smooth $k$-algebra of relative dimension $n$, 
such that the $A$-module 
$T=Der_k(A)$ is free and admits a base $\{\btau_i\}$ consisting of commuting 
vector fields, and let $E$ be a free $A$-module of rank $m$, with a base 
$\{\phi_{\alpha}\}$. 
We shall call the set ${\mathfrak g}=\{\btau_i;\;\phi_{\alpha}\}\subset T\oplus E$ 
a {\it frame} of $(T,E)$. 
Consider the commutative $A$-superalgebra 
$\Lambda E=\bigoplus_{i=0}^m\Lambda_A^i(E)$ where the parity of 
$\Lambda^i_A(E)$ is equal to the parity of $i$. Each frame $\mathfrak g$ as above 
gives rise  to a vertex superalgebroid as follows.

Let $T_{\Lambda E}=Der_k(\Lambda E)$.
We extend the fields $\btau_i$ to derivations $\tau_i$ of the whole superalgebra 
$\Lambda E$ by the rule 
$$
\tau_i(a)=\btau_i(a),\qquad
\tau_i(\sum a_{\alpha}\phi_\alpha)=\sum \btau_i(a_{\alpha})\phi_{\alpha}
.
$$
(Note that this extension  depends on a choice 
of a base $\{\phi_{\alpha}\}$ of the module $E$, though the whole construction will be independent of the basis by the work of~\cite{chiral3,msv}.) The fields $\{\tau_i\}$ 
form a $\Lambda E$-base of the even part $T_{\Lambda E}^{ev}$. 

We define the odd vector fields $\psi_{\alpha}\in T_{\Lambda E}^{odd}$ by 
$$
\psi_{\alpha}(\sum a_{\nu}\phi_{\nu})=a_{\alpha},\qquad\psi_{\alpha}(a)=0.
$$
These fields form a $\Lambda E$-base of $T_{\Lambda E}^{odd}$. 

Let $\{\omega_i;\;\rho_{\alpha}\}$ be the dual base of the module 
of $1$-superforms $\Omega_{\Lambda E}=
Hom^{ev}_{\Lambda E}(T_{\Lambda E},\Lambda E)$, 
defined by 
$$
\langle\tau_i,\omega_j\rangle=\delta_{ij},\qquad
\langle\psi_{\alpha},\rho_{\beta}\rangle=\delta_{\alpha\beta},\qquad
\langle\tau_i,\rho_{\alpha}\rangle=
\langle\psi_{\alpha},\omega_i\rangle=0.
$$ 
Maps $c$ and $\gamma$ are completely determined by setting them to zero on the elements of each frame. They extend uniquely to give a vertex superalgebroid structure (see \cite[section 3.4]{chiral3}).
Thus  we have constructed a vertex superalgebroid $(\Lambda(E),T_{\Lambda E},\Omega_{\Lambda E},\partial,
\gamma,\langle\ ,\ \rangle,c)$.

Now we can construct the vertex superalgebroids whose vertex envelopes give our two main examples: the vertex algebras of polyvector fields and chiral de Rham forms.

Let us work in the analytic category. Let $k$ be the field of complex numbers $\mathbb C$ and let $A$ be the algebra of analytic functions on $\mathbb C^n$. Let $E$ be the module of vector fields $T$ or the module of $1$-forms $\Omega^1_{A/k}$ over $\mathbb C^n$;  
its exterior algebras 
are the algebra of polyvector fields and the de Rham algebra of differential forms 
$\Omega^{\bullet}_{A/k}$ over $\mathbb C^n$. These choices of $E$ produce vertex algebroids according to the construction above.
Their vertex envelopes are the required topological vertex algebras called in \cite{chiral3} the chiral de Rham complex and chiral vector field complex over $\mathbb C^n$. 
In \cite{msv} a sheafification construction was presented which produces from the local objects we have just described a sheaf of topological vertex algebras on a (smooth, analytic or algebraic) manifold. 

By our main result, Theorem \ref{maintheorem}, the local objects have the structure of $G_\infty$-algebras.
Therefore there is a natural task to investigate whether they pass to a global structure on the manifold using the glueing procedure of \cite{msv}.

In conformal weight zero we have 
a global Gerstenhaber structure in the polyvector field case; 
in the chiral de Rham complex there is no bracket and 
the structure is just the associative algebra of differential forms.
We always have a global operations $x_{(-1)}y$, $x_{(0)}y$ and $G_{(0)}$ 
giving the Lian--Zuckerman product and bracket, and 
our inductive construction will give an extension to a global $G_\infty$-structure on the chiral de Rham sheaf as long as 
the fields corresponding to $L$, $G$, $Q$ and $J$ extend to global sections of the appropriate sheaf of vertex algebras. 
By  \cite[Theorem 4.2]{msv} this will occur if the first Chern class $c_1$ of the manifold vanishes.

\begin{cor}
Let $M$ be a $C^\infty$ Calabi--Yau manifold. Then there is a structure of a $G_\infty$-algebra on the vector space underlying the chiral de Rham complex or the complex of chiral polyvector fields on $M$ which is locally the $G_\infty$ structure constructed above.
\end{cor}

One can see that the above homotopy structure descends to the
cohomology of the chiral sheaves making it a homotopy
Gerstenhaber algebra too.

\subsection{An $L_{\infty }$ algebra related to the Courant bracket} 

Recall that Dorfman
\cite{dorfman} and Courant \cite{courant} defined the following brackets on the space of sections of the bundle $TM+T^{\ast }M$ on a manifold $M$: 
\begin{eqnarray*}
\left[ X+\psi ,Y+\nu \right] _{D} &=&
\left[ X,Y\right] +\left( L_{X}\nu -\iota _{Y}d\psi \right)\\ \left[ X+\psi
,Y+\nu \right] _{C} &=&{\textstyle\frac12}\left( \left[ X+\psi ,Y+\nu \right] _{D}%
-\left[ Y+\nu ,X+\psi \right] _{D}\right) 
\end{eqnarray*}%
where we follow the notation of \cite{gualtieri}: $X,Y$ are vector fields and $%
\psi $ and $\nu $ are $1$-forms. The Courant bracket is a key ingredient
in recent developments in generalized complex geometry \cite{gualtieri,hitchin}. This bracket is skewsymmetric, but the Jacobi identity
holds only modulo the image of the de Rham differential. It poses a natural
problem of finding an $L_{\infty }$ algebra structure $(l_1,l_2,l_3,\ldots)$ whose $l_{1}$ is the de Rham
differential, and whose $l_{2}$ is related to the Courant bracket. Such an
algebra does not have to be unique of course. A solution to the problem
was given in~\cite{rw-courant}; the $L_{\infty }$-algebra constructed
there has $l_{m}=0$ for all $m\geq4$.

In the presence of supersymmetry in the language of conformal field
theories, the following refinement of the above problem makes sense. Suppose
the de Rham operator splits as a supercommutator of two operators.
Then either of these operators can be taken as $l_{1}$ in the above
problem. 

We address this new refined problem in the context of the chiral de
Rham complex, which posesses the required supersymmetry in the case of
manifolds whose first Chern class is zero. In this case, as explained above,
the chiral de Rham complex is a (sheaf of) topological vertex algebras
whose operator $L_{(0)}$ splits as a supercommutator, $L_{(0)}=[Q_{(0)},G_{(0)}]$. On the other hand in the chiral de Rham complex the
operator $L_{(0)}$ restricted to the functions on the manifold is just the classical de Rham
differential. 

To make the connection with the Courant bracket recall from~\cite[Section 3.10]{msv} that the chiral de Rham complex posesses a
filtration, whose associated graded object {\bf Gr} is identified. 
{\bf Gr} is also called the quasiclassical limit of the chiral de Rham complex.
Now $TM+T^{\ast }M$ is a summand in {\bf Gr} of the component
of conformal weight one, and a beautiful observation by Bressler~\cite{bressler-p1} states that the quasiclassical limit of the operation $x_{(0)}y$ is
the Courant bracket when restricted to this summand. 

This suggests that
the graded piece {\bf Gr}$_1$ of conformal weight one is an $L_{\infty }$ algebra
which solves our problem. The only modification necessary to the theory we have developed for the Lian--Zuckerman conjecture is the following. The bracket, defined above by means of
the composite operation $(G_{(0)}x)_{(0)}y$, must be replaced by the operation $x_{(0)}y$ alone. This is already a skewsymmetric operation on {\bf Gr}$_1$, and trivially satisfies the Jacobi identity in cohomology since the complex is exact in non-zero conformal weights. 
Our argument, or that of~\cite{shlie}, therefore produces an $L_\infty$ structure on {\bf Gr}$_1$ starting 
with $l_1x=Q_{(0)}x$ as before, and $l_{2}(x,y)=x_{(0)}y$.
By \cite{msv,bressler-p1} we have: 

\begin{prop}
\begin{enumerate}
\item The piece of {\bf Gr}$_1$ of conformal weight one of {\bf Gr} is a module over $\Pi T^{*}M$ generated by 
$$TM,\qquad T^{*}M,\qquad \Pi TM,\qquad \Pi T^*M_{(-2)}\b1$$
where $\b{1}$ is the vacuum vector and $\Pi $ is the change of parity functor.
\item The restriction of the operation $x_{(0)}y$ to the submodule
generated by $TM\oplus T^{*}M$ is equal to the Courant bracket.
\item The operation $l_2(x,y)=x_{(0)}y$ satisfies the Jacobi identity in cohomology (relative to boundary $d=l_1$) and an inductive construction as in Theorem \ref{main} produces an $L_{\infty }$-algebra structure on {\bf Gr}$_1$ which extends $l_2$ and hence the Courant bracket.
\end{enumerate}\end{prop}

This construction of an $L_{\infty }$-algebra related to the
Courant bracket is natural since it is derived from of the chiral de Rham
complex, a part of the geometric structure of a manifold.

\end{document}